\documentclass[aop,secthm,dvips]{arximspdf}

\doi{10.1214/009117906000000980}
\volume{35}
\issue{3}
\pubyear{2007}
\firstpage{915}
\lastpage{934}

\makeatletter

\newcommand{\nn}{\nonumber}
\newcommand{\ff}{\iffalse}

\newproclaim{remark}[thm]{Remark}
\newtheorem{lemma}[thm]{Lemma}
\newtheorem{corr}[thm]{Corollary}

\newcommand{\real}{\mathbb{R}}
\makeatother

\begin{document}
\begin{frontmatter}

\title{On the continuity of local times of Borel right Markov processes}
\runtitle{Continuity of local times}

\begin{aug}
\author[A]{\fnms{Nathalie} \snm{Eisenbaum}\ead[label=e1]{nae@ccr.jussieu.fr}} and
\author[B]{\fnms{Haya} \snm{Kaspi}\corref{}\ead[label=e2]{iehaya@tx.technion.ac.il}}
\runauthor{N. Eisenbaum and H. Kaspi}
\affiliation{Universit\'{e} Paris VI-CNRS and Technion}
\address[A]{Laboratoire de Probabilit\'es\\
            Universit\'e Paris VI---CNRS\\
            4 place Jussieu\\
            75252 Paris Cedex 05\\
            France\\
            \printead{e1} } %adresu isvedimo komanda gale!
\address[B]{Industrial Engineering\\
            \quad and Management\\
            Technion\\
            Technion City\\
            Haifa 32000\\
            Israel\\
            \printead{e2}}
\end{aug}

% HISTORY:
\received{\smonth{12} \syear{2005}}
\revised{\smonth{10} \syear{2006}}

% ABSTRACT
\begin{abstract}
The problem of finding a necessary
and sufficient condition for the continuity of the local times for
a general Markov process  is still open.   Barlow and Hawkes  have
completely treated the case of the L\'evy processes, and Marcus and
Rosen  have solved the case of the strongly symmetric Markov
processes. We treat here the continuity of the local times of Borel
right processes. Our approach unifies that of Barlow and Hawkes and of Marcus
and Rosen, by using an associated Gaussian process, that appears as a limit in a
CLT involving the local time process.
\end{abstract}

% KEYWORDS

\end{frontmatter}

%s1 ###
\section{Introduction}

Let $X=(\Omega,\mathcal{F},\mathcal{F}_t,X_t,\theta_t,P_x; x\in
E)$ be a Borel right process, having a reference measure $m$, with
all states communicating and regular for themselves. Under these
assumptions, a local time $L^x_t$ exists at each point, unique up
to a multiplicative constant. Let $u^{\alpha}(x,y)$ be the
potential densities with respect to $m$, and normalize the local
times (choose the multiplicative constant), so that for some (and
all) $\alpha,$
%e1 ###
\begin{equation}\label{e1}
E_x\int_0^\infty e^{-\alpha t}\,dL^y_t=u^{\alpha}(x,y),
\end{equation}
for all $x,y\in E$, where $E_x$ is the expectation with respect to
$P_x$. The question, under what conditions there exists  a
version of
 $(L^x_t)_{x\in E, t>0}$ so that $(x,t)\to L^x_t(\omega)$ is almost surely
 continuous, has occupied many researchers in the field for many years.
 Although, as we shall describe below, there are
 some very important special cases where this problem has been
 solved, the problem, for general Borel right processes, is still
 unresolved.

To put this problem in context, we would like to start by
highlighting some of the most important existing results in this
field.
 The first to address this problem was Trotter who in
\cite{Tr} proved that when $X$ is the Brownian motion on the
real line, it has a local time at all points and (normalized as
above) there is a version of $(x,t)\to L^x_t$ that is almost
surely jointly continuous. In \cite{GK} Getoor and Kesten
have treated the problem for standard Markov processes that have a
reference measure. They have established a sufficient condition and
a necessary condition for the above joint continuity, but with a
gap between the necessary and the sufficient conditions.

 Bass and Khoshnevisan in \cite{BK}, Barlow in
\cite{B1,B2} and Barlow and Hawkes in
\cite{BH} have treated the case of L\'{e}vy processes taking
real values that have local times at all points and for which all
points communicate. In \cite{B2} the
 necessary and sufficient conditions for the
existence of  an almost surely, jointly continuous version of the
local time $(x,t)\to L^x_t(\omega)$ were found. Since this solution is in
many ways the starting point of our approach, we shall describe
it here (or rather Bertoin's \cite{Be} ``translation'' of it).

Let
$h(a,b)=E_a(L^a_{T_b})=E_0(L^0_{T_{b-a}})=h(0,b-a)\mathop{=}\limits^{\mathrm{def}}h(b-a)$.
Then one can show that $h(x)=h(-x)$ and that $d^2(a,b)=h(b-a)$
defines a distance on $\real$ that is equivalent to the Euclidean
distance. Let $m(y)=|\{x\dvtx h(x)<y\}|$, where $|A|$ is the Lebesgue
measure of the Borel set $A\subset\real$. Barlow's necessary and
sufficient condition  for the continuity of the local time is the
following ``majorizing measure'' condition:
%e2 ###
\begin{equation}\label{e2}
\int_{0+}^{\bolds\cdot}\sqrt{\ln\frac{1}{m(\varepsilon)}}\,d\varepsilon<\infty.
\end{equation}
In Barlow's paper \cite{B2} the  condition is stated in terms of
the inverse of $m$,  the monotone rearrangement of $h$. Let
$\bar{h}(x)=\inf\{y\dvtx m(y)>x\}$, then $(x,t)\to L^x_t$ has a
continuous version, iff
%e3 ###
\begin{equation}\label{e3}
I(\bar{h}) =\int_{0+}^{\bolds\cdot}\frac{\bar{h}(x)}{x(\ln
x)^{1/2}}\, dx<\infty.
\end{equation}
It is easily seen that (\ref{e2}) and (\ref{e3}) are equivalent, but as was
noticed by Barlow and Hawkes \cite{BH}, (\ref{e3}) is reminiscent of
Fernique's  \cite{F} and Dudley's \cite{Du} necessary and
sufficient condition on the covariance function of a stationary
Gaussian process $(\phi_x)$ to have a continuous version. This
Gaussian process is precisely described in \cite{E}.

In a series of papers during the 1990s \cite{MR1,MR2,MR3,MR4,MR5}
and in their recent book \cite{MR6}, Marcus and Rosen study
sample path properties of the local time process of strongly
symmetric Markov processes. Under symmetry, the potential
densities are symmetric and positive definite. Therefore, there
exists a centered Gaussian process $(\phi_x)_{x\in E}$ such that
$\langle\phi_x\phi_y\rangle =u(x,y)$, where $u(x,y)$ is the 0-potential
density when the process is transient and $\langle\phi_x\phi_y\rangle=u^{\alpha}(x,y)$ for
$\alpha>0$ when the process is recurrent.  Now and in the
sequel, $\langle \cdot \rangle$ denotes the expectation with respect to the Gaussian measure.
The main tool for their study is the celebrated Dynkin isomorphism
theorem (DIT) \cite{D1,D2}, which states when $X$ is transient  for any measurable function $F$
on $\real^{E}$,
%e4 ###
\begin{equation}\label{e4}
E_{a,b}\biggl\langle F\biggl(L^{\bolds\cdot}_{\zeta}+\frac{\phi_{\bolds\cdot}^2}{2}\biggr) \biggr\rangle
= \biggl\langle\frac{\phi_a\phi_b}{\langle\phi_a\phi_b\rangle}
F\biggl(\frac{\phi_{\bolds\cdot}^2}{2}\biggr)\biggr\rangle,
\end{equation}
where $\zeta$ is the life time of the  Markov process
$X$  and  $E_{a,b}$ is the law of $X$ born at $a$ and killed at its
last exit from the point $b$. Note that when $X$ is recurrent, the above identity is available
for $X$ killed at an independent exponential time with parameter $\alpha$. One should notice that the
right-hand
side of (\ref{e4}) is stated in terms of the Gaussian process only.

Defining a distance by
$d^2(x,y)=\langle(\phi_x-\phi_y)^2\rangle =u(x,x)-2u(x,y)+u(y,y)$, they have
used the DIT to show that $(x,t)\to L^x_t$ has a jointly
continuous version (in the distance $d$), iff the Gaussian process
$(\phi_x)$ has a continuous version in that distance. The latter
happens iff for every compact set $K$, in the metric $d$, there
exists a probability measure $\mu$ on $\tilde{\mathcal{K}}$, the
$\sigma$-algebra on $K$ generated by the $d$-open sets, so that
%e5 ###
\begin{equation}\label{e5}
\lim_{\delta\to 0}\sup_{x\in
  K}\int_0^{\delta}\sqrt{\ln\frac{1}{\mu(B(x,\varepsilon))}}\, d\varepsilon=0,
\end{equation}
where $B(x,\varepsilon)$ is a $d$-ball of radius $\varepsilon$
around $x$. When $(x,y)\to u(x,y)$ is jointly continuous those
conditions translate to a condition for
the joint continuity (in the original distance) of $(x,t)\to
L^x_t$. As in the case treated by Barlow, the condition for the
joint continuity of the local time is identical to the condition
for the continuity of a Gaussian process.

 Extending these results beyond the symmetric and L\'{e}vy cases, and
understanding the intriguing connection between the conditions for
the continuity of local times of Markov processes and those of
Gaussian processes is the objective of this paper.

We shall work under the following assumptions:
\begin{longlist}[(A4)]
\item[(A1)] All points of $E$ are regular for themselves.

\item[(A2)] All points of $E$ communicate.

\item[(A3)] The process is recurrent.

\item[(A4)] There exists a Borel right dual process.
\end{longlist}
The recurrence property will simplify our
arguments considerably, but it is not a very serious assumption.
Indeed, by an argument due to Le Jan (see \cite {DM}, Chapter~XII),
if $X$ is transient, one can always ``revive'' it in such a
way that it becomes recurrent, still keeping properties that will
be used below like duality or symmetry if the original process
was symmetric. Since the continuity and other fine properties of
local times are local, and recurrence is a long time behavior
property of the process, it has nothing to do with local
properties, and therefore, using Le Jan's construction, we can
extend the results to the transient case. With that in mind, and
assuming that $X$ is recurrent, let $m$ be the unique invariant
distribution for $(P_t).$ (A1) and (A2) imply that $m$ is actually
a reference measure. Thus, the potential densities
$u^{\alpha}(x,y)$ exist. From general theory (see \cite{Fi}), we
know that a dual process $\hat{X}$ exists. That is, there exists a
Markov process $\hat{X}$ whose potential is given by
$\hat{U}^{\alpha}f(x)=\int_E m(dy)f(y)u^{\alpha}(y,x).$ Since $X$
is recurrent, so is $\hat{X}$, and since $X$ has a local time at
each $x\in E$, so does $\hat{X}$. In general $\hat{X}$ is not a
strong Markov process. It is only a moderate Markov process,
namely, it satisfies the strong Markov property only at
($\mathcal{F}_t)$ predictable stopping times. Our fourth
assumption (A4) and the only serious one (beyond those needed to
define the problem properly) is that $\hat{X}$ is actually a
Borel right process as well, or that at least it satisfies the
strong Markov property at the hitting times $T_x \ \mbox{of all}\
x\in E$. Note that the L\'{e}vy processes treated by Barlow satisfy
this assumption (with $\hat{X}=-X$), and the symmetric processes
studied by Marcus and Rosen satisfy it with $\hat{X}=X.$

 To state our main results, we shall need  some
additional notation.  Let $0$ be a preassigned state in $E$ and
$T_0$ be its hitting time. By recurrence, $T_0<\infty$ $P_x$ a.s.
for every $x\in E$. Let $u_{T_0}(x,y)$ be the potential densities
of the process $X^{T_0}$, where
\[
 X^{T_0}_t=\cases{   X_t,    & \quad if $t<T_0$, \cr
                   \Delta, &\quad otherwise,}
\]
where $\Delta$ is a cemetery state. $X^{T_0}$ is the process
killed
 at its hitting time of the state~0. We shall show  that
$u_{T_0}(x,y)+u_{T_0}(y,x)$ is both symmetric and positive
definite. Thus, there exists a centered Gaussian process $(\phi_x)_{x\in E}$, such that
$\langle\phi_x\phi_y\rangle=u_{T_0}(x,y)+u_{T_0}(y,x).$ Using this, we now
define the distance $d$  with which we shall work:
%e6 ###
\begin{equation}\label{e6}
\quad d^2(x,y)=u_{T_0}(x,x)-u_{T_0}(x,y)-u_{T_0}(y,x)+u_{T_0}(y,y)= \langle(\phi_x-\phi_y)^2\rangle.
\end{equation}

 Our first result gives a sufficient condition for the continuity
 of the local time process.
\begin{thm}\label{t1.1}
 If for every compact set $K$, in the $d$ metric, there exists a
probability measure $\mu$ on the Borel sets of $K$ defined with
the $d$-distance, so that
\[
\lim_{\delta\to 0}\sup_{x\in
  K}\int_0^{\delta}\sqrt{\ln\frac{1}{\mu(B(x,v))}}\,dv=0,
  \]
  where $B(x,\varepsilon)$ is a ball of $d$-radius $\varepsilon$ around
 $x$,
   then $(t,x)\to L^x_t(\omega)$ has a jointly $d$-continuous version. If further
  $(x,y)\to  u^{\alpha}(x,y)$ is jointly continuous, then a.s. $(x,t)\rightarrow L^x_t$ is
  continuous in $d$ and the original distances. Finally, for a compact set $K$, set
 \[
 \eta_K(\delta)=\sup_{z\in
K}\int_0^{\delta}\sqrt{\ln\frac{1}{\mu(B(z,v))}}\,dv.
\]
%Then, if $\lim_{\varepsilon\rightarrow0}\frac{\eta_K(\varepsilon)}{\varepsilon} \rightarrow \infty$
 There is a positive
  constant C such that
  \[
  \limsup_{\delta\rightarrow 0}\sup_{a,b\in K
 d(a,b)<\delta}\sup_{s\le t}\frac{|L^a_s-L^b_s|}{\eta_K(d(a,b))}\le
 C\biggl(\sup_{x\in K}L^x_t\biggr)^{1/2}.
 \]
\end{thm}

   Since the sufficient condition of Theorem \ref{t1.1} is actually a necessary and sufficient
condition for the continuity of the Gaussian process $\phi$ (see \cite{LT}),
Theorem~\ref{t1.1} contains the  following relation:

\noindent{\sl If $(\phi_x)_{x\in E}$ has a  continuous version for the distance $d$\emph{,}
then $(L^x_t, x\in E, t\geq 0)$ has a jointly continuous version for the distance $d$.}

 Our next two theorems deal with a central limit
theorem in $C(K)$, the space of continuous functions on a compact set
$K$ contained in $E$. We believe that this theorem provides the missing link
for the converse of the above relation.

 Let $\tau^a(s)=\inf\{t>0\dvtx L^a_t>s\}$. Then
$\tau^a(s)$ is a process with stationary independent increments.
In particular, $\tau^0(s)$ is a process with stationary
independent increments, and $L^{\bolds\cdot}_{\tau^0(s)}$ is a process
with stationary independent increments taking values in function
space. We note further that, for any $s$, $L^{\bolds\cdot}_{\tau^0(s)}$
has an infinitely divisible law, and therefore,
$Y_n(\cdot)=\frac{L^{\bolds\cdot}_{\tau^0(n)}-n}{\sqrt{n}}$ is an
infinitely divisible random variable, taking values in the space
of functions. We refer the reader to \cite{AG} and \cite{L} for more on infinitely
divisible processes taking values in Banach spaces.

\begin{thm}\label{t1.2}
 If $x\to Y_1(x)$ is a continuous function in the $d$ distance, and
  the majorizing measure condition of Theorem \emph{\ref{t1.1}} holds,
  then for each compact set $K$ in the $d$ metric,
$(Y_n(x))_{x\in K}$ converges weakly in $C(K)$ to a centered Gaussian
process $(\phi_x)_{x\in K}$ satisfying
$\langle\phi_x\phi_y\rangle=u_{T_0}(x,y)+u_{T_0}(y,x)$.
\end{thm}

 The characterization of continuous Gaussian processes as those
 for which the covariance distance satisfies the above majorizing
 measure condition yields the following theorem as a corollary.

\begin{thm}\label{t1.3}
Assume that $u^{\alpha}(x,y)$ are continuous,
 then the following are equivalent:
\begin{longlist}[2.]
\item[1.] $(x,t)\rightarrow L^x_t$ is jointly continuous and the above CLT holds.
\item[2.] The above majorizing measure condition holds.
\end{longlist}
\end{thm}

\begin{remark}\label{r1.4}
 We have not been able to show
that the continuity of the local time process alone is a sufficient
condition for the majorizing measure condition to hold.
 However, in view of all
existing results, we conjecture that this is really the case.
Theorem \ref{t1.3} allows one to replace the proof of sufficiency with a
proof that the continuity of $x\rightarrow Y_n(x)$ implies its
tightness in $C(K)$.
\end{remark}

 Our paper is organized as follows. In Section \ref{s2} we
prove some preliminary results on the metric $d(x,y)$ defined in
(\ref{e6}). Those will be our main tool for proving our results. Section
\ref{s3} is devoted to the proof of Theorem \ref{t1.1} and in Section \ref{s4} we shall
prove Theorem \ref{t1.2}  and Theorem \ref{t1.3} as its corollary. We shall also
recall there from \cite{EK} that in the symmetric case the
tightness that is needed for the CLT follows easily from the
continuity of the associated Gaussian process and the results of
\cite{EKMRS}.

%s2 ###
\section{Notation and preliminary results}\label{s2}

 We adopt the basic notation of Blumenthal and Getoor \cite{BG2}. We let
$X,\hat{X}$ be two recurrent Borel right Markov processes in classical
duality. As can be easily shown, under (A1)--(A3), the unique
invariant measure $m$ for this process is also a reference
measure.   Let $u^{\alpha}(x,y)$ be the corresponding potential densities,
 $U^{\alpha}f(x)=\int u^{\alpha}(x,y)f(y)m(dy)$ and
$\hat{U}^{\alpha}f(y)=\int u^{\alpha}(x,y)f(x)m(dx)$. Hence,
$u^{\alpha}(x,y)$ is the potential density of the process $X$
starting at $x$ and $u^{\alpha}(x,y)$ is the potential density of
the process $\hat{X}$ starting at $y$. We will assume from now on
that the processes have local times at each point (enough to assume
that one of them has a local time at each point, the other will have
it as a result), and that the local times are normalized so that
\[
 u^{\alpha}(x,y)=E_x\int e^{-\alpha t}\,dL^y_t
 \]
and similarly for the dual process,
\[
u^{\alpha}(x,y)=\hat{E}_y\int e^{\alpha t}\,d\hat{L}^x_t.
\]

For every state $x\in E$, let $T_x=\inf\{t>0\dvtx X_t=x\}$, we shall use
the notation $T_x$ for the dual process as well. Denote by
$u_{T_x}(a,b)$ the potential densities of the process killed at
$T_x$. The two  resulting processes are again in duality with respect
to $m(dy)$. By recurrence, $u_{T_x}(a,b)$ is finite and is equal
to the increasing limit of $u^{\alpha}_{T_x}(a,b)$ as $\alpha\to
0$. Let $\nu^x$ be the excursion measure from $x$ and similarly
for the dual process, denote it by $\hat{\nu}^x$. All excursions
from a point end at this point.
\begin{lemma}\label{l1}
Let $x,y$ be two points in $E$. Then
$u_{T_x}(y,y)=u_{T_y}(x,x)$.
\end{lemma}
\begin{pf}
Recall that
$u_{T_x}(y,y)=E_y(L^y_{T_x})=\lim_{\alpha \to
0}E_y\int_0^{T_x}e^{-\alpha t}\,dL^y_t$.

Therefore,
\[
\frac{u_{T_x}(y,y)}{u_{T_y}(x,x)}=\lim_{\alpha \to
0}\frac{E_y\int_0^{T_x}e^{-\alpha
t}\,dL^y_t}{E_x\int_0^{T_y}e^{-\alpha t}\,dL^x_t}.
\]
Now,
\begin{eqnarray*}
u^{\alpha}(x,x)&=&E_x\int_0^\infty e^{-\alpha
t}\,dL^x_t
 \\
&=&E_x\int_0^{T_y}e^{-\alpha t}\,dL^x_t + E_x(e^{-\alpha
T_y})E_y(e^{-\alpha T_x})u^{\alpha}(x,x).
\end{eqnarray*}
Hence,
%e7 ###
\begin{equation}\label{e7}
E_x\int_0^{T_y}e^{-\alpha t}\,dL^x_t=u^{\alpha}(x,x)\bigl(1-E_x(e^{-\alpha
T_y})E_y(e^{-\alpha T_x})\bigr)
\end{equation}
 and similarly,
%e8 ###
 \begin{equation}\label{e8}
 E_y\int_0^{T_x}e^{-\alpha
t}\,dL^y_t=u^{\alpha}(y,y)\bigl(1-E_y(e^{-\alpha T_x})E_x(e^{-\alpha
T_y})\bigr).
\end{equation}
Our result will follow if we can show that
\[
\lim_{\alpha\to 0}\frac{u^{\alpha}(x,x)}{u^{\alpha}(y,y)}=1.
\]
But,
\[
\frac{u^{\alpha}(x,x)}{u^{\alpha}(y,y)}=\frac{u^{\alpha}(x,x)}{u^{\alpha}(y,x)}
\frac{u^{\alpha}(y,x)}{u^{\alpha}(y,y)}=\frac{\hat{E}_x(e^{-\alpha
T_y})}{E_y(e^{-\alpha T_x})}.
\]
Since $X$ is recurrent, so is the dual $\hat{X}$ and thus,
\[
\lim_{\alpha\to 0}\frac{\hat{E}_x(e^{-\alpha
T_y})}{E_y(e^{-\alpha T_x})}=
\frac{\hat{P}_x(T_y<\infty)}{P_y(T_x<\infty)}=1.
\]\upqed
\end{pf}
With this result at hand we now have the following:
\begin{lemma}\label{l2.2}
For every $x$,$y$ in $E$,
 \[
 \nu^0\bigl((L^x-L^y)^2\bigr)=2\bigl(u_{T_0}(x,x)-u_{T_0}(x,y)-u_{T_0}(y,x)+u_{T_0}(y,y)\bigr).
\]
\end{lemma}
\begin{pf}
Let $(\theta_t)$ be the usual shift operators on the
state space so that $X_s(\theta_t\omega)=X_{t+s}(\omega)$ and
$\hat{\theta}_t$ defined similarly for the dual process:
%e9 ###
 \begin{equation}\label{e9}
\nu^0(L^xL^y)=\nu^0\biggl(\int_0^{T_0}L^y_{T_0}(\theta_
t)\,dL^x_t+\int_0^{T_0}L^x_{T_0}(\theta_t)\,dL^y_t\biggr).
\end{equation}
 By the Markov property that $\nu^0$ satisfies, this is equal to
\[
\nu^0\biggl(\int_0^{T_0}E_x(L^y_{T_0})\,dL^x_t+\int_0^{T_0}E_y(L^x_{T_0})\,dL^y_t\biggr)
\]
and hence, to
\[
\nu^0\bigl(1_{\{T_x<T_0\}}u_{T_0}(x,x)u_{T_0}(x,y)+1_{\{T_y<T_0\}}u_{T_0}(y,y)u_{T_0}(y,x)\bigr).
\]
But
%e10 ###
\begin{equation}\label{e10}
\nu^0(T_x<T_0)=\frac{1}{E_0(L^0_{T_x})}=\frac{1}{u_{T_x}(0,0)}=\frac{1}
{u_{T_0}(x,x)},
\end{equation}
 where the last equality follows from Lemma \ref{l1}. Inserting this into (\ref{e9})
 yields
%e13 ###
%e12 ###
%e11 ###
 \begin{eqnarray}\label{e11}
\nu^0(L^xL^y)&=&\frac{u_{T_0}(x,x)u_{T_0}(x,y)}{u_{T_0}(x,x)}\nn
\\
&&{}+
\frac{u_{T_0}(y,y)u_{T_0}(y,x)}{u_{T_0}(y,y)}
\\
&=&u_{T_0}(x,y)+u_{T_0}(y,x).\nn
\end{eqnarray}
\upqed
\end{pf}
\begin{corr}\label{c3}
$(u_{T_0}(x,y)+u_{T_0}(y,x), x,y \in E\times E)$ is
symmetric, positive definite.
\end{corr}
\begin{pf}
Symmetry is obvious. Let $(a_1,\ldots,a_n)$ be a vector in
 $\real^n$, then
\begin{eqnarray*}
&& \sum_{i=1}^n\sum_{j=1}^n
 a_ia_j\bigl(u_{T_0}(x_i,x_j)+u_{T_0}(x_j,x_i)\bigr)
 \\
 &&\qquad =\nu^0\Biggl(\Biggl(\sum_{i=1}^n
 a_iL^{x_i}\Biggr)^2\Biggr)\ge 0.
\end{eqnarray*}\upqed
\end{pf}
 We now define
 \[
 d^2(x,y)=u_{T_0}(x,x)-u_{T_0}(x,y)-u_{T_0}(y,x)+u_{T_0}(y,y).
 \]
The above results prove that $d(x,y)$ defines a pseudo distance,
and that there is a centered Gaussian process $(\phi_x)$ such that
%e14 ###
\begin{equation}\label{e12}
\langle\phi_x,\phi_y\rangle=u_{T_0}(x,y)+u_{T_0}(y,x).
\end{equation}
\begin{lemma}\label{l4}
Set $h(x,y)=E_x(L^x_{T_y})$, then
$d^2(x,y)=h(x,y)$.
\end{lemma}
\begin{pf}
\begin{eqnarray*}
d^2(x,y)&=&u_{T_0}(x,x)-u_{T_0}(x,y)-u_{T_0}(y,x)+u_{T_0}(y,y)\\
&=& h(x,0)+h(y,0)-P_x(T_y<T_0)h(y,0)-P_y(T_x<T_0)h(x,0)\\
&=& h(x,0)P_y(T_x>T_0)+h(y,0)P_x(T_y>T_0)\\
&=& u_{T_0}(x,x)\frac{u_{T_x}(y,0)}{u_{T_x}(0,0)}+u_{T_0}(y,y)
\frac{u_{T_y}(x,0)}{u_{T_y}(0,0)}.
\end{eqnarray*}
But by Lemma \ref{l1}, $u_{T_0}(x,x)=u_{T_x}(0,0)$ and
$u_{T_0}(y,y)=u_{T_y}(0,0)$, and the last term is equal to
\begin{eqnarray*}
&& u_{T_x}(y,0)+u_{T_y}(x,0)
\\
&&\qquad =\hat{E}_0(L^y_{T_x})+\hat{E}_0(L^x_{T_y})
\\
&&\qquad =\hat{P}_0(T_y<T_x)\hat{E}_y(L^y_{T_x})+\hat{P}_0(T_x<T_y)\hat{E}_x(L^x_{T_y})
\\
&&\qquad =\hat{E}_y(L^y_{T_x})=\hat{h}(y,x),
\end{eqnarray*}
where the one before last equality follows from Lemma \ref{l1} applied
to the dual process. We now notice that since $d^2(x,y)$ is
symmetric with respect to $x$ and $y$, and with respect to the
dual objects, it follows that
$\hat{h}(y,x)=\hat{h}(x,y)=h(y,x)=h(x,y)$, and our result
follows.
\end{pf}
\begin{remark}\label{r5}
It follows from the
above result that $d$ is a real distance on $E$. Indeed, for $x\ne
y$, $P_x(T_y > 0 )= 1$ and since $x$ is regular for itself, this
implies that $E_x(L^x_{T_y}) > 0$ and hence, that $d(x,y) > 0$. If
the potential densities $u_{T_0}(x,y)$ are jointly continuous,
continuity in the topology generated by this metric implies
continuity in the original metric on $E$.
\end{remark}

%s3 ###
\section{Sufficiency of the majorizing measure condition}\label{s3}

Thanks to the results of the previous section, the proof of sufficiency is
very close to that of Bertoin's (\cite{Be}, pages 144--150).
\begin{lemma}\label{l3.1}
For
$a\in E$, set $\tau^a_t=\inf\{s\dvtx  L^a_s>t\}$. Then for every $a,b\in
E$,
%e15 ###
\begin{equation}\label{e13}
P\{\exists s\le \tau^b_y \dvtx L^b_s-L^a_s>x\}\le
\exp\biggl(-\frac{x^2}{4yh(a,b)}\biggr).
\end{equation}
\end{lemma}
\begin{pf}
By the Markov property, we may start our process at
$b$. $t\to L^a_{\tau^b_t}$ is a subordinator (which may have a jump
at $0$ if our process does not start at $b$). $L^a_{\tau^b_t}$ stays
at $0$ and performs its first jump at time $L^b_{T_a}$. It has
therefore no drift. $L^b_{T_a}$ under $P_b$ has an exponential
distribution with expectation $E_b(L^b_{T_a})=h(a,b)$. Next let
$R=\inf\{t>T_a\dvtx X_t=b\}$ and note that $R=\tau^b_{L^b_{T_a}}$. Since $L^a_{T_a}=0$,
$L^a_R$ has again an exponential distribution with expectation
$E_a(L^a_{T_b})=h(a,b)$. Hence, the L\'{e}vy measure of $L^a_{\tau^b_t}$
is equal to $\frac{1}{h^2(a,b)}\exp(-\frac{1}{h(a,b)}x)$, its L\'{e}vy
exponent $\Psi(\lambda)$ is  equal to $\frac{\lambda}{\lambda
h(a,b)+1}$, and $\exp(-\lambda
L^a_{\tau^b_s}+\frac{s\lambda}{\lambda h(a,b)+1})$ is a martingale.
By the optional sampling theorem applied to $T\wedge y$ where
$T=\inf\{s\dvtx  s-L^a_{\tau^b_s}> x\}$, we can show that
\[
P\{T\le y\}\le \exp \biggl(-\lambda x+\lambda y\biggl(1-\frac{1}{\lambda
h(a,b)+1}\biggr)\biggr).
\]
 Taking now $ \lambda=\frac{x}{2yh(a,b)}$ gives us the required upper
 bound.
 \end{pf}

 Define now
%e16 ###
 \begin{equation}\label{e14}
 Y_a(q)(t)=q\wedge L^a_t.
 \end{equation}
Then for every $q>0, a,b\in E$, if $q\wedge L^b_t-q\wedge L^a_t>
x$ at some time $t\ge 0$, then the time when this occurs is
bounded by $\tau^b_q$. Hence,
\begin{eqnarray*}
&&\{\exists s\dvtx |Y_a(q)(s)-Y_b(q)(s)|>x\}
\\
&&\qquad = \{\exists s\le \tau^b_q\dvtx
L^b_s-L^a_s>x\}\cup\{\exists s\le \tau^a_q: L^a_s-L^b_s>x\}.
\end{eqnarray*}
Therefore,
%e17 ###
\begin{equation}\label{e15}
P\{|Y_a(q) - Y_b(q)|_u > x\}\le 2\exp\biggl(-\frac{x^2}{4qh(a,b)}\biggr),
\end{equation}
where $|\ |_u$ is the uniform bound with respect to time. It now
follows that, for every $c>0$,
\begin{eqnarray*}
&& E\bigl(\exp\bigl(|Y_b(q) -Y_a(q)|^2_u/c\bigr)-1\bigr)
\\
&&\qquad =\frac{1}{c}\int_0^\infty\exp\biggl(\frac{x}{c}\biggr)
P\{|Y_b(q)-Y_a(q)|^2_u>x\}\,dx
\\
&&\qquad \le \frac{2}{c}\int_0^\infty
\exp\biggl(\frac{x}{c}\biggr)\exp\biggl(\frac{x}{4qh(a,b)}\biggr)\,dx = 2\biggl(\frac{c}{4qh(a,b)}-1\biggr)^{-1}.
\end{eqnarray*}
 Taking now $c>12qh(a,b)$, we get
\[
E\bigl(\exp\bigl(|Y_b(q)-Y_a(q)|^2_u\bigr)-1\bigr)<1.
\]
 In the language of Ledoux and Talagrand (page 298 in \cite{LT}),
 \[
\|Y_b(q)-Y_a(q)\|_{\psi}\le \tilde{d}(a,b),
\]
where $\tilde{d}^2(a,b)=12qh(a,b)$ and the Young function
 $\psi(x)=\exp(x^2)-1$. We shall fix now $q>1$ and  abuse the
 notation by denoting $ Y_a(q)(t)$ by $Y_a(t)$.
\begin{pf*}{Proof of Theorem \ref{t1.1}}
\textit{Step} 1.
By Theorem 11.14 of \cite{LT}, if
 for a compact set $K$ there exists a probability measure $\mu$ on $(K,\tilde{d}$) such that for
 \[
\lim_{\eta\to 0} \sup_{x\in
 K}\int_0^{\eta}\sqrt{\ln{\frac{1}{\mu(B(x,\varepsilon))}}}\, d\varepsilon=0,
 \]
where $B(x,\varepsilon)$ is a ball of radius $\varepsilon$ in the distance
 $ \tilde{d}$ around $x$, then  for each $x\in E$, $(a\to Y_a)_{a\in K}$
  has a $P_x$ almost surely
 continuous version with respect to the distance $\tilde{d}$.
 That is, there is a process$(\tilde{Y}_a)_{a\in K}$ with
continuous sample paths in the distance $\tilde{d}$ (and therefore
$d$), so that for every $a\in K$, $Y_a=\tilde{Y}_a$  $P_x$ almost
surely.

 Let $(\tilde{Y}_a)_{a\in K}$ be that version and define
 $\bar{Y}_{*}(t)=\sup\{\tilde{Y}_a(t)\dvtx a\in K\}$ and
 $\Theta(q)=\inf\{t\dvtx \bar{Y}_{*}(t)=q\}$. Since
 $a\rightarrow\tilde{Y}_a$
 is continuous,
 $\bar{Y}_{*}(t)=\sup\{\tilde{Y}_r(t): r \in \Gamma\}$, where $\Gamma$ is a countable
 dense set in $K$ with respect to the distance $\tilde{d}$. Thus,
 $\bar{Y}_{*}(t)$ is a nondecreasing adapted process. Hence, $\Theta
(q)$ is an $(\mathcal{F}_t)$ stopping time. By the Blumenthal 0--1
law, $\{\Theta(q)=0\}$ is a probability $0$ or $1$ event for every
$P_z$.

\textit{Step} 2. \textit{$\Theta(q)>0$ for $q$ large enough}.
We shall show that, for $q$ large enough,
$\Theta(q)>0,$ $P_x$ a.s. and that $\Theta(q)\rightarrow \infty $
as $q\rightarrow\infty.$ To do that, we shall use Proposition 1
of \cite{H}. Indeed, let
\[
f(q,\omega, a,b)=
|\tilde{Y}_a(\tau(1))-\tilde{Y}_b(\tau(1))|
\]
 and
\[
\tilde{f}(q,\omega, a,b)=\frac{f(q,\omega,
a,b)}{\tilde{d}(q,a,b)},
\]
 where
 \[
\tau(1)=\inf\{s\dvtx \tilde{Y}_0(s)\ge 1\},
\]
 where $0$ is a preassigned
state, that we assume is in $K$. Then
 \begin{eqnarray*}
 && P_x\{\exp(\tilde{f}^2(q,\omega,a,b))-1>\alpha\}
 \\
 &&\qquad =P_x\{\exp(\tilde{f}^2(q,\omega,a,b))>1+\alpha\}
 \\
 &&\qquad = P_x\{f^2(q,\omega,a,b)> 12\ln(1+\alpha)qd^2(a,b)\}
 \\
 &&\qquad =P_x\bigl\{f(q,\omega,a,b)>\sqrt{12q\ln(1+\alpha)}d(a,b)\bigr\}
 \\
 &&\qquad \le P_x\biggl\{\sup_{s\le\tau^0_1}|L^a_s\wedge q-L^b_s\wedge q|>\sqrt{12q\ln(1+\alpha)}d(a,b)\biggr\}
 \\
 &&\qquad \le  2\exp\biggl(-\frac{12q\ln(1+\alpha)d^2(a,b)}{4qh(a,b)}\biggr)
 \\
 &&\qquad \le 2\exp\biggl(-\frac{12q\ln(1+\alpha)}{4q}\biggr)
 \\
 &&\qquad =2\biggl(\frac{1}{1+\alpha}\biggr)^3,
 \end{eqnarray*}
 where the first inequality  follows from (\ref{e15}) and the third from
 Lemma \ref{l2.2}. It follows that, for all $q$,
%e18 ###
 \begin{equation}\label{eq16}
 E_x(\exp(\tilde{f}^2(q,\omega,a,b)))\le 3 .
 \end{equation}
Define now
\[
C(q)=\int_K\int_K
 \exp(\tilde{f}^2(q,\omega,a,b))\mu(da)\mu(db).
 \]
Then $E_x(C(q))\le 3$, and therefore, $C(q)<\infty$, $P_x \mbox{ a.s.}$
Since $a\rightarrow \tilde{Y}_a(\tau(1))$ is continuous in the
$\tilde{d}(=\tilde{d}(q))$ distance, we can use Heinkel's formula
to deduce that, for all $(a,b)\in K$,
\[
|\tilde{Y}_a(\tau(1))-\tilde{Y}_b(\tau(1))|\le
20\sup_{z\in
K}\int_0^{\sqrt{12q}d(a,b)/2}\biggl(\ln\biggl(\frac{C(q)}{\mu^2(\tilde{B}(z,u))}\biggr)\biggr)^{1/2}\,du,
\]
where, as before, $\tilde{B}(x,v)$ is a ball of radius $v$ in the
$\tilde{d}$ distance. This after the change of variable
$v=\frac{u}{\sqrt{3q}}$ is equal to
\[
20\sqrt{3q}\sup_{z\in
K}\int_0^{d(a,b)}\biggl(\ln\biggl(\frac{C(q)}{\mu^2(B(z,v/2))}\biggr)\biggr)^{1/2}\,dv,
\]
where
$B(z,v)$ is a ball of radius $v$ in the $d$ distance. This is
bounded above by
\[
40\sqrt{3q}\sup_{z\in
K}\int_0^{d(a,b)}\biggl(\ln\biggl(\frac{C(q)}{\mu^2(B(z,v))}\biggr)\biggr)^{1/2}\,dv,
\]
which is easily shown to be bounded by
\[
40\sqrt{3q}(\ln C(q))^{1/2}d(a,b)+\sqrt{2}\sup_{z\in
K}\int_0^{d(a,b)}\biggl(\ln\biggl(\frac{1}{\mu(B(z,v))}\biggr)\biggr)^{1/2}\,dv.
\]

For $\delta>0$, define
%e19 ###
\begin{equation}\label{e17}
\eta(\delta)=\sup_{z\in
K}\int_0^{\delta}\biggl(\ln\biggl(\frac{1}{\mu(B(z,v))}\biggr)\biggr)^{1/2}\,dv.
\end{equation}
By our assumptions, $\eta(D)<\infty$, where $D$ is the diameter of
the compact set $K$, and $\lim_{\delta\rightarrow
0}\eta(\delta)=0$. Returning to our computation,
\[
\tilde{Y}_a(\tau(1))\le \tilde{Y}_0(\tau(1))+c\sqrt{q}\bigl((\ln
C(q))^{1/2}D+\sqrt{2}\eta(D)\bigr),
\]
where $c$ is a constant and $D$ is
the diameter of $K$ in the $d$ distance. We now recall that $P_x$
a.s.  $\tilde{Y}_0(\tau(1))=1$, and so, on $\{\Theta(q)\le
\tau(1)\}$,
\[
q\le 1+c\sqrt{q}\bigl((\ln C(q))^{1/2}D+\sqrt{2}\eta(D)\bigr).
\]
Using, as in \cite{B1}, the fact that $y^2\le A+By$ implies that
$y^2\le 2A+B^2$, we see that, on $\{\Theta(q)\le \tau(1)\}$,
\[
q \le 2+\bigl(c((\ln C(q))^{1/2}D+\sqrt{2}\eta(D))\bigr)^2.
\]
Thus,
\begin{eqnarray*}
P_x\{\Theta(q)\le \tau(1)\}&\le & P_x\bigl\{\bigl(c((\ln
C(q))^{1/2}D+\sqrt{2}\eta(D))\bigr)^2\ge q-2\bigr\}
\\
&\le & P_x\biggl\{2((\ln
C(q))^{1/2}D)^2+4\eta^2(D)\ge \frac{q-2}{c^2}\biggr\}
\\
&=& P_x\biggl\{\ln
C(q)\ge \frac{q-2}{2c^2D^2}-\frac{2\eta^2(D)}{D^2}\biggr\}
\\
&=&
P_x\biggl\{C(q)\ge
\exp\biggl(\frac{q-2}{2c^2D^2}-\frac{2\eta^2(D)}{D^2}\biggr)\biggr\}
\\
&\le& 3\exp\biggl(-\frac{q-2}{2c^2D^2}+\frac{2\eta^2(D)}{D^2}\biggr),
\end{eqnarray*}
where the last inequality follows from the fact that $E_x(C(q))\le
3$, this last term is smaller than $1$ for $q$ large enough. It
now follows that $P_x\{\Theta(q)=0\}\le P_x\{\Theta(q)\le
\tau(1)\}<1$, so that $P_x\{\Theta(q)=0\}=0$ for $q$ large enough.\vadjust{\goodbreak}

\textit{Step} 3. $\lim_{n\rightarrow\infty}\Theta(n)$.
\noindent Repeating  the above computation with
$q(n)=3n$ and $\tau(n)=\inf\{t>0\dvtx  Y_0(t)>n\}$, instead of $q = 1$ and $\tau(1)$, we similarly obtain
\begin{eqnarray*}
P_x\{\Theta(q(n))\le \tau(n)\} &\le &
3\exp\biggl(-\frac{3n-2n}{2c^2D^2}+\frac{2\eta^2(D)}{D^2}\biggr)
\\
&\le & A\exp\biggl(-\frac{n}{B}\biggr),
\end{eqnarray*}
where A and B are constants. It now follows from the fact that
$\tau(n)\rightarrow \infty$  as $ n\rightarrow \infty$, and the
Borel--Cantelli lemma that $\Theta(3n)\rightarrow \infty,$ $P_x$
a.s., and so, the local time has a jointly continuous version
(with respect to $P_x$) in $\real_+\times K$ for all compact K,
and thus a jointly continuous version.

\textit{Step} 4. \textit{Jointly continuous potential
densities}.
 To prove the last assertion of Theorem \ref{t1.1}, we
need to show that when $u^{\alpha}(x,y)$ is jointly continuous,
then a.s. $(t,x)\rightarrow L^x(t,\omega)$ is jointly continuous.
First note that by Remark 3.4.4 of \cite{MR6} $u^{\beta}(x,y)$ is
jointly continuous for any $\beta>0$, and thus, it can be easily
shown that $u_{T_0}(x,y)$ is jointly continuous as well. Hence, the
continuity in the $d$ distance implies continuity in the original
distance [Actually, as was noted by Bertoin (\cite{Be}, page 147),
this will make the distances equivalent in the sense that their
induced topologies will be the same.]
%Indeed,  $x$ close to $y$ in
%the $d$ distance, implies that $E_x(L^x_{T_y})$, is small, which
%implies that $T_y$ is small when one starts from $x$. By the right
%continuity of $t\rightarrow X_t$ in the original distance, $y$
%must be close to $x$ in that distance).
We fix now a $q$ and prove
that, $P_x$ almost surely, $\tilde{Y}_a(q)(t)=L^a_t$
simultaneously for all $a\in K$ and $t<\Theta(q)$. Specifically,
for each $a\in K,$ $P_x$ almost surely, $L^a_t=\tilde{Y}_a(q)(t)$
for all $t<\Theta(q)$. By Fubini's theorem and the occupation time
density formula \cite{BG1}, for each continuous $f$ with support
in $K$,
\[
 \int_0^tf(X_s)\,ds=\int_Kf(a)\tilde{Y}_a(q)(t)m(da)
 \]
for all $t<\Theta(q)$ $P_x$ almost surely. Letting $f$ range over
a countable dense family in $C(K)$ (the continuous functions with
support in $K$), we see that the identity above holds for all
bounded measurable functions with support in $K$. In particular,
$P_x$ almost surely for all $a\in K$, $\varepsilon>0$ and
$t<\Theta(q)$,
\[
\int_0^tf_{\varepsilon,a}(X_s)\,ds=\int_Kf_{\varepsilon,a}(x)\tilde{Y}_x(q)(t)m(dx),
\]
where $f_{\varepsilon,a}$ is an approximating delta function that
defines the local time at $a\in K$, (see Theorem 3.6.3 and the
discussion preceding it in Marcus and Rosen's recent
book \cite{MR6}, that can be easily adapted to the nonsymmetric
situation). Since $x\to \tilde{Y}_x(q)(t)$ is continuous, the
right-hand side converges to $\tilde{Y}_a(q)(t)$ as $\varepsilon\to 0$,
and $\lim_{\varepsilon \to 0}\int_0^tf_{\varepsilon,a}(X_s)\,ds=L^a_t$
uniformly in $[0,\Theta(q))$ by Theorem 3.6.3 of~\cite{MR6}. Thus,
$\tilde{Y}_a(q)(t)=L^a_t$, $P_x$ almost surely for all $a\in K$
and $ t<\Theta(q)$, so that $(a,t) \to L^a_t$ is $P_x$ a.s.
continuous on $K\times[0,\Theta (q)))$. Since we have seen that
$P_x$ a.s. $\Theta(3n)\rightarrow\infty$, as $
n\rightarrow\infty$, it follows that $(a,t) \to L^a_t$ is $P_x$
a.s. continuous on $K\times[0,\infty)$ and therefore, in
$E\times[0,\infty)$. Since this is true for every $x\in E$, it
follows that a.s. $(a,t) \to L^a_t$ is continuous.

\textit{Step} 5. \textit{Modulous of continuity}.
 To get the modulous result, we follow \cite{B1},
with Heinkel's inequality replacing the Gracia, Rodemich, Rumsey
inequality that appears there. We now return to Heinkel's
inequality with $L^x_t$, which we now know is continuous in the
$d$ distance. Let $K$ be a compact set in the $d$
metric. For a fixed $t >0$,  let $q$ be such that $\sup_{x\in K} L^x_t(\omega)\le q\le
2\sup_{x\in K}L^x_t(\omega)$. Note that we have
\[
\liminf_{\delta\downarrow 0}{1\over \delta}\eta_K(\delta) \geq \lim_{\delta \downarrow 0}
\biggl(\ln {1\over \mu(B(z,\delta))} \biggr)^{1/2}
\]
for every $z$ in $K$. Hence, unless $\mu$ charges all the points of $K$ and $K$ is
finite,
we obtain
\[
\lim_{\delta\downarrow  0}{1\over \delta}\eta_K(\delta) = + \infty
\]
[we have used the fact that for every infinite countable subset $A$ of
$K$,\break
$\sum_{z \in A} \mu(\{z\}) \leq 1$]. Of course, in the case when $K$ is finite,
the question of the modulous of continuity is meaningless. Hence, we can choose
 $\delta(\omega)>0 $  small enough such that
 $(\ln(
C(q,K,\omega)))^{1/2}\delta(\omega)<\sqrt{2}\eta_K(\delta(\omega))$.
If  $\varepsilon_t(\omega)$ is chosen to be smaller than
$\delta(\omega)$, then, as in the computation of Step 2, with
$(L^a_t)_{a\in K}$ replacing $(\tilde{Y}_a(\tau(1))_{a\in K})$
that appears there, one can show that for $ a,b \in K $,
\[
|L^a_s-L^b_s|\le  C\biggl(\biggl(\sup_{x\in K}L^x_t\biggr)^{1/2}\eta_K(d(a,b)\biggr)
\]
for
all $0\le s\le t$  and all $ a,b\in K$ such that
$d(a,b)<\varepsilon_t(\omega)$. $C$ is a constant, which by our
computations is smaller than 80, and is by no means the best
possible (see \cite{B1,B2,BK}  in the L\'{e}vy case and
\cite{MR1} in the symmetric case). We shall not pursue this issue
any further here.
\end{pf*}

%s4 ###
\section{The central limit theorem for local times}\label{s4}

Trying to understand the true reason why in all
existing results the conditions for the joint continuity of the local
time are identical to
those for the continuity of Gaussian processes, one is led to seek the
explanation in a suitable CLT. Indeed, let  $\tau_t^0$ be the
inverse of the local time at $0$, then $(L^x_{\tau^0_t})_{x\in E}$,
is a process with stationary independent increments with values
that
are functions on $E_0$. By Lemma \ref{l2.2} and its proof,
$E(L^x_{\tau^0_n})=n$, and  $E_0((L^x_{\tau^0_n})^2)=2u_{T_0}(x,x)n$.
It follows that $\frac{L^x_{\tau^0_n}-n}{\sqrt{n}}$ converges in
distribution to a centered normal random variable with variance
$2u_{T_0}(x,x)$. The following lemma will show that the process
$(\frac{L^x_{\tau^0_n}-n}{\sqrt{n}})_{x\in E}$ converges in finite-dimensional
distributions to a Gaussian process with covariance
$u_{T_0}(x,y)+u_{T_0}(y,x)$.
\begin{lemma}\label{l4.1}
\[
E_0\biggl(\frac{(L^x_{\tau^0_n}-n)(L^y_{\tau^0_n}-n)}{n}\biggr)=
u_{T_0}(x,y)+u_{T_0}(y,x).
\]
\end{lemma}
\begin{pf}
 Denote by $G$ the set of left endpoints of
excursions from
$\{t\dvtx X_t=0\}$. Then we have
%e20 ###
\begin{equation}\label{e18}
E_0(L^x_{\tau^0_n})=E_0\Biggl(\sum_{s\in
G,s\le\tau^0_n}L^x_{T_0}\circ\theta_s\Biggr)=E_0\int_0^{\tau^0_n}\nu^0(L^x)\,dL^0_t=n,
\end{equation}
where the second equality follows from excursion theory
(compensating the sum of jumps), and the third by a change of
variable $s=L^0_t$, and the fact that $\nu^0(L^x)=1$ for all $x\in
E$, which we have shown in the proof of Lemma \ref{l2.2}. Similarly,
%e23 ###
%e22 ###
%e21 ###
\begin{eqnarray}\label{e19}
E_0(L^x_{\tau^0_n}L^y_{\tau^0_n}) &=& E_0\sum_{s\in
G,s\le\tau^0_n}L^x_{T_0}\circ\theta_s\sum_{t\in
G,t\le\tau^0_n}L^y_{T_0}\circ\theta_t\nn
\\
&=& E_0\sum_{s\in G, s\le
\tau^0_n}(L^x_{T_0}L^y_{T_0})\circ\theta_s
\\
&&{}+E_0\sum_{s\in G, s\le
\tau^0_n}L^x_{T_0}\circ\theta_s \sum_{t\in G, t\le \tau^0_n, t\neq
s}L^y_{T_0}\circ \theta_t.\nn
\end{eqnarray}
Using excursion theory as above, the first sum of (19) is equal to
%e24 ###
\begin{equation}\label{e20}
E_0\int_0^{\tau^0_n}\nu^0(L^xL^y)\,dL^0_t = n\nu^0(L^xL^y)
 =
n\bigl(u_{T_0}(x,y)+u_{T_0}(y,x)\bigr),
\end{equation}
 where the last equality follows as in the proof of Lemma \ref{l2.2}. Using
 excursion theory again, the
 second term of (\ref{e19}) is composed of two
 sums. The first is equal to
 \[
E_0\int_0^{\tau^0_n}\nu^0(L^x)E_0(L^y_{\tau^0_n-t})\,dL^0_t
\]
 and the second is identical to the above with $x$ and $y$ interchanged. Since
  \mbox{$\nu^0(L^x)=1$},  this integral is
 equal to $ E_0\int_0^nE_0(L^y_{\tau^0_{n-u}})\,du=
 E_0\int_0^nE_0(L^y_{\tau^0_u})\,du=\frac{n^2}{2} $. Since the value of
 this term is independent of $x$ and $y$, it remains the same when
 interchanging $x$ and $y$. Thus,
 \begin{eqnarray*}
&&E_0 \biggl(\frac{(L^x_{\tau^0_n}-n)(L^y_{\tau^0_n}-n)}{n}\biggr)
\\
&&\qquad =\frac{n(u_{T_0}(x,y)+
u_{T_0}(y,x))+n^2-n^2-n^2+n^2}{n}
\\
&&\qquad =u_{T_0}(x,y)+ u_{T_0}(y,x).
\end{eqnarray*}\upqed
\end{pf}

\begin{pf*}{Proof of Theorem \ref{t1.2}}
All one needs to prove is tightness in $C(K)$ which,
since $ L^0_{\tau^0_n}-n=0,$
amounts to showing that for every $\eta>0, \varepsilon>0, \exists
\delta>0$, \mbox{$\exists n_0\in \mathbf{N}$}, so that for all $n\ge n_0$,
\[
P\biggl\{\sup_{a,b\in K,d(a,b)<\delta}\frac{|L^a_{\tau^0_n}-L^b_{\tau^0_n}|}{\sqrt{n}}>
\eta\biggr\} <\varepsilon.
\]
To prove this, we shall use here
again Proposition 1.1 of \cite{H}, but in view of the
computations in the last section, with some steps abridged. We
shall split the proof into a few steps.

\textit{Step \emph{1}. Definition of objects appearing in Heinkel\emph{'}s
inequality}:
%e26 ###
%e25 ###
\begin{eqnarray}
Y_n(q,a,b)&=&\sup_{s\le\tau^0_n}|q\wedge L^a_s-q\wedge L^b_s|,\label{e21}
\\
\tilde{Y}_n(q,a,b)&=&\frac{Y_n(q,a,b)}{\rho(a,b)},\label{e22}
\end{eqnarray}
where $\rho(a,b)=\sqrt{12q}d(a,b)$. Then as in the computations
preceding (16), for all $y\in E,$
%e27 ###
\begin{equation}\label{e23}
P_y\{\exp(\tilde{Y}_n^2(q,a,b))-1>\alpha\} \le
2\biggl(\frac{1}{1+\alpha}\biggr)^3.
\end{equation}
Hence, for all $n$, $E_y(\exp(\tilde{Y}_n^2(q,a,b)))\le 3$.

Define now
%e28 ###
\begin{equation}\label{e24}
C(n,q)= \int_K\int_K \exp(\tilde{Y}_n^2(q,a,b))\mu(da)\mu(db)
\end{equation}
Then $E_y(C(n,q))\le3$ and therefore, $C(n,q)<\infty,$ $P_y$
a.s. for all $y\in E$.

\textit{Step \emph{2}. Application of Heinkel\emph{'}s inequality}.
It follows from Proposition 1.1 of~\cite{H} that
\[
|q\wedge L^a_{\tau^0_n}-q\wedge L^b_{\tau^0_n}|\le 20 \sup_{z\in
K}\int_0^{\sqrt{12q}d(a,b)/2}\biggl(\ln\biggl(\frac{C(n,q)}{\mu^2(B(z,u))}\biggr)\biggr)^{1/2}\,du.
\]
Following again the same arguments as in the previous section, the
right-hand side of the above inequality is bounded above by
\[
40\sqrt{3q}\biggl((\ln C(n,q))^{1/2}d(a,b)+\sqrt{2}\sup_{z\in
K}\int_0^{d(a,b)}\biggl(\ln\biggl(\frac{1}{\mu(B(z,v))}\biggr)\biggr)^{1/2}\,dv\biggr).
\]
Recall the definition of $\eta(\delta)$ from (\ref{e17}). By our
assumption, $\eta(\delta)\rightarrow 0$ as $\delta \rightarrow
0$.\vadjust{\goodbreak}

 \textit{Step \emph{3}.}
With the above result at hand, we now take
$q=n+\lambda\sqrt{n}$ to obtain
\begin{eqnarray*}
&& \frac{|(n+\lambda\sqrt{n})\wedge
L^a_{\tau^0_n}-(n+\lambda\sqrt{n})\wedge
L^b_{\tau^0_n}|}{\sqrt{n}}
\\
&&\qquad \le \frac{1}{\sqrt{n}}40\sqrt{3}\bigl(n+
\lambda\sqrt{n}\bigr)^{1/2}
\bigl[\bigl(\ln\bigl(C\bigl(n,n+\lambda\sqrt{n}\bigr)\bigr)\bigr)^{1/2}d(a,b)+\sqrt{2}\eta(d(a,b))\bigr]
\\
&&\qquad =40\sqrt{3}\biggl(1+\frac{\lambda}{\sqrt{n}}\biggr)^{1/2} \bigl[\bigl(\ln\bigl(
C\bigl(n,n+\lambda\sqrt{n}\bigr)\bigr)\bigr)^{1/2}d(a,b)+\sqrt{2}\eta(d(a,b))\bigr].
\end{eqnarray*}

Returning now to the proof of tightness,
%e31 ###
%e30 ###
%e29 ###
\begin{eqnarray}
&& P \biggl(\sup_{d(a,b)<a,b\in K,
\delta}\frac{|L^a_{\tau^0_n}-L^b_{\tau^0_n}|}{\sqrt{n}}>\eta\biggr)\label{e25}
\\
 &&\qquad \le  P\biggl(\sup_{a,b\in K, d(a,b)<\delta}\frac{|(n+\lambda\sqrt{n})\wedge
L^a_{\tau^0_n}-(n+\lambda\sqrt{n})\wedge
L^b_{\tau^0_n}|}{\sqrt{n}}>\eta\biggr)\label{e26}
\\
 &&\qquad\quad {}+ P\biggl(\sup_{x\in
 K}L^x_{\tau^0_n}>n+\lambda\sqrt{n}\biggr).\label{e27}
\end{eqnarray}

 \textit{Step} 4.
 Starting with (\ref{e27}) and using the above inequality
with $a=x, b=0$, and recalling that $L^0_{\tau^0_n}=n$, we get
\[
\bigl(n+\lambda\sqrt{n}\bigr)\wedge L^x_{\tau^0_n}\le
n+40\sqrt{3}\bigl(n+\lambda\sqrt{n}\bigr)^{1/2}\bigl[\bigl(\ln\bigl(
C\bigl(n,n+\lambda\sqrt{n}\bigr)\bigr)\bigr)^{1/2}D+
\sqrt{2}\eta(D)\bigr],
\]
where $D$ is the diameter of $K$ with respect to the distance $d$
which, by our assumption, is finite as is $\eta(D)$. Now on
$\{\sup_{x\in K}L^x_{\tau^0_n}>n+\lambda\sqrt{n}\}$,
\[
n+\lambda\sqrt{n}\le
n+40\sqrt{3}\bigl(n+\lambda\sqrt{n}\bigr)^{1/2}\bigl[\bigl(\ln\bigl(
C\bigl(n,n+\lambda\sqrt{n}\bigr)\bigr)\bigr)^{1/2}D+ \sqrt{2}\eta(D)\bigr],
\]
so that
\[
\lambda\sqrt{n}\le 40\sqrt{3}\bigl(n+\lambda\sqrt{n}\bigr)^{1/2}\bigl[\bigl(\ln\bigl( C\bigl(n,n+\lambda
\sqrt{n}\bigr)\bigr)\bigr)^{1/2}D+\sqrt{2}\eta(D)\bigr],
\]
which is equivalent  to
\[
\frac{\lambda}{(1+\lambda/\sqrt{n})^{1/2}}\le 40\sqrt{3}\bigl[
\bigl(\ln\bigl( C\bigl(n,n+\lambda \sqrt{n}\bigr)\bigr)\bigr)^{1/2}D+\sqrt{2}
\eta(D)\bigr].
\]
Thus,
\begin{eqnarray*}
&& P\biggl\{\sup_{x\in K}L^x_{\tau^0_n}>n+\lambda\sqrt{n}\biggr\}
\\
&&\qquad \le P\biggl\{40\sqrt{3}\bigl[\bigl(\ln\bigl( C\bigl(n,n+\lambda
\sqrt{n}\bigr)\bigr)\bigr)^{1/2}D+\sqrt{2}\eta(D)\bigr]\ge
\frac{\lambda}{(1+\lambda/\sqrt{n})^{1/2}}\biggr\}.
\end{eqnarray*}
This last probability is equal to
%e32 ###
\begin{equation}\label{e28}
P\biggl\{\ln\bigl(C\bigl(n,n+\lambda\sqrt{n}\bigr)\bigr)\ge
\biggl(\frac{\lambda}{(1+\lambda/\sqrt{n})^{1/2}40\sqrt{3}D}-\frac{\sqrt{2}
\eta(D)}{D}\biggr)^2\biggr\}.
\end{equation}
 For $\lambda$ big enough, so
that the first term on the right-hand side of the inequality of (\ref{e28}) is
larger than $3$ times the second, this is smaller or equal to
\begin{eqnarray*}
&& P\biggl\{\ln\bigl(C\bigl(n,n+\lambda\sqrt{n}\bigr)\bigr)\ge
\frac{\lambda^2}{3(1+\lambda/\sqrt{n})4800D^2}\biggr\}
\\
&&\qquad \le
3\exp\biggl(-\frac{\lambda^2C}{D^2(1+ \lambda/\sqrt{n})}\biggr)\le
A\exp\biggl(-\frac{\lambda^2}{1+\lambda}B\biggr),
\end{eqnarray*}
where the first inequality follows from the fact that we have
shown that $E(C(n, q))\le 3$ for all $q$ and $n$, and $A,B,C$ are
some constants.

 We now choose $\lambda$ big enough to satisfy all the above
inequalities and make this last bound smaller than
$\varepsilon/2$. Note that this $\lambda$ is chosen independently
of $n$. With this $\lambda$, we return to (\ref{e26}).

 \textit{Step} 5.
\begin{eqnarray*}
&& P\biggl\{\sup_{d(a,b)<\delta}\frac{|(n+\lambda\sqrt{n})\wedge
L^a_{\tau^0_n}-(n+\lambda\sqrt{n})\wedge
L^b_{\tau^0_n}|}{\sqrt{n}}>\eta\biggr\}
\\
&&\qquad \le P\biggl\{40\sqrt{3}\biggl(1+\frac{\lambda}{\sqrt{n}}\biggr)^{1/2} \bigl[\bigl(\ln\bigl(
C\bigl(n,n+\lambda\sqrt{n}\bigr)\bigr)\bigr)^{1/2}\delta+\sqrt{2}\eta(\delta)\bigr]>\eta\biggr\}
\\
&&\qquad =P\biggl\{\bigl(\ln\bigl(
C\bigl(n,n+\lambda\sqrt{n}\bigr)\bigr)\bigr)^{1/2}\delta+\sqrt{2}\eta(\delta)>\frac{\eta}{40
\sqrt{3}(1+ \lambda/\sqrt{n})^{1/2}}\biggr\}.
\end{eqnarray*}
Choose now $\delta^*$ small  enough so that
$\sqrt{2}\eta(\delta)<\frac{\eta}{80\sqrt{3}(1+\lambda)^{1/2}}$
for all $\delta\le\delta^*$. For $\delta\le\delta^*$, the above
probability is smaller or equal to
\[
P\biggl\{\bigl(\ln\bigl(C\bigl(n,n+\lambda\sqrt{n}\bigr)\bigr)\bigr)^{1/2}\delta>
\frac{\eta}{80\sqrt{3}(1+\lambda/\sqrt{n})^{1/2}}\biggr\},
\]
which is equal to
\begin{eqnarray*}
&&P\biggl\{\bigl(\ln\bigl( C\bigl(n,n+\lambda\sqrt{n}\bigr)\bigr)\bigr)>
\frac{\eta^2c}{(1+\lambda/\sqrt{n})\delta^2}\biggr\}
\\
&&\qquad \le P\biggl\{\bigl(\ln\bigl( C\bigl(n,n+\lambda\sqrt{n}\bigr)\bigr)\bigr)>
\frac{\eta^2c}{(1+\lambda)\delta^2}\biggr\}
\\
&&\qquad \le 3\exp\biggl(-\frac{\eta^2c}{(1+\lambda)\delta^2}\biggr),
\end{eqnarray*}
where $c$ is a constant. Note that this bound is independent of $n$ and one can choose
$\delta$ small enough to satisfy all the above inequalities and
make it smaller than  $\ \varepsilon/2$,  which proves the desired
tightness.
\end{pf*}
\setcounter{thm}{2}
\begin{remark}\label{r4.3}
Both $K$ and $C(K)$ are defined with respect
to the metric $d$. If the potential densities are jointly
continuous, this will imply a corresponding CLT with respect to
the original metric.
\end{remark}
\begin{remark}\label{r4.4}
The following was done in
\cite{EK} when the process $X$ is symmetric; we present it here
again for the sake of completeness.   It has been shown in~\cite{EKMRS} that
%e33 ###
\begin{equation}\label{e29}
\bigl(L^x_{\tau_n}+\tfrac{1}{2}\phi_x^2; x\in E\bigr)\mathop{=}^{\mathrm{law}}\bigl(\tfrac{1}{2}\bigl(\phi_x+
\sqrt{n}\bigr)^2; x\in E\bigr).
\end{equation}
Subtracting $n$ from both sides and dividing them by $\sqrt{n}$,
\[
\biggl(\frac{L^x_{\tau_n}-n}{\sqrt{n}}+\frac{\phi_x^2}{2\sqrt{n}};
x\in E\biggr)\mathop{=}^{\mathrm{law}}
\biggl(\frac{\phi^2_x}{2\sqrt{n}}+\sqrt{2}\phi_x; x\in E\biggr)
\]
and our tightness in $C(K)$ follows directly from the tightness of the Gaussian law.
 See also \cite{AMZ} for tightness in the symmetric case using the DIT directly.
\end{remark}
\begin{pf*}{Proof of Theorem \ref{t1.3}}
Under the assumption that
$u^{\alpha}(x,y)$ are continuous, Theorems \ref{t1.1} and \ref{t1.2} show that
the existence of a majorizing measure is sufficient for both the
continuity and the tightness in $C(K)$ of the local time process
both with respect to the metric $d$ and then, by the above
continuity, with respect to the original distance on $E$. The
necessity follows from the characterization of Gaussian processes
that are continuous in the metric $d$, as those for which the
majorizing measure conditions are satisfied for each compact set
$K$.
\end{pf*}

\printaddresses

\end{document}